\documentclass{amsart}

\usepackage[applemac]{inputenc}   

\newtheorem{theorem}{Theorem}[section]

\newtheorem{corollary}[theorem]{Corollary}

\theoremstyle{definition}

\theoremstyle{remark}

\numberwithin{equation}{section}

\begin{document}

\title{Dynamical Zeta Functions and Kummer Congruences}

\author{J. Arias de Reyna}
\address{Department of Mathematical Analysis, 
University of Seville,
Seville,
Spain}
\curraddr{Facultad de Matemáticas, Universidad de Sevilla, Apdo.~1160, 
41080-Sevilla, Spain}
\email{arias@us.es}
\thanks{This research was supported  by  Grant BFM2000-0514.}


\subjclass{Primary 11B68, 37C30; Secondary 11B37, 37B10}

\date{June 26, 2003 }


\keywords{Kummer congruences, Bernoulli numbers, Euler numbers, 
integer sequences, zeta function}

\begin{abstract}
We establish a connection between the coefficients of Artin-Mazur
zeta-functions and Kummer congruences. 
\par
This allows to settle positively the question of the existence of a 
map $T\colon X\to X$ such that the number of fixed 
points of $T^n$ are $|E_{2n}|$, where $E_{2n}$ are the Euler numbers. 
Also we solve a problem of Gabcke related to the coefficients of Riemann-Siegel 
formula.
\end{abstract}

\maketitle

\def\N{\mathbf{N}}
\def\Z{\mathbf{Z}}
\def\Fix{\operatorname{Fix}}
\def\sgn{\operatorname{sgn}}
\def\Tr{\operatorname{Tr}}

\section*{Introduction}

In this paper we establish a connection between two important topics 
the Artin-Mazur zeta function and Kummer's congruences.  Some 
connection between Kummer's congruences and periodic points are 
pointed in the paper by Everest, van der Poorten, Puri and Ward 
\cite{E}.

Inspired by the Hasse-Weil zeta function of an algebraic variety over a 
finite field, Artin and Mazur \cite{AM} defined  the Artin-Mazur zeta 
function for an arbitrary map $T\colon X\to X$ of a topological space  $X$:
$$Z(T;x):=\exp\left(\sum_{n=1}^\infty\frac{\Fix T^n}{n} 
x^n\right).$$
Where $\Fix T^n$ is the number of isolated fixed points of 
$T^n$. 

Manning \cite{M}  proved the rationality of the Artin-Mazur zeta 
function for diffeomorphisms of a smooth compact manifold satisfying 
Smale axiom $A$. 

Following \cite{P}, call a sequence $a=(a_n)_{n\ge1}$ of 
non-negative integers \emph{ realizable} if there is a set $X$ and a map 
$T\colon X\to X$ such that $a_n$ is the number of fixed points of 
$T^n$. 

We must notice that in \cite{P} it is proved that if $(a_{n})$ is 
realizable, then there exists a compact space $X$ and a  homeomorphism 
$T\colon X\to X$, such that $a_{n}=\Fix T^n$.

Puri and Ward \cite{PW} proved that a sequence of non-negative 
integers 
 $(a_{n})_{n\ge1}$  is realizable if and only if 
$\sum_{d\mid n}\mu(n/d)a_d$
is non negative and divisible by $n$ for all $n\ge1$.
Here $\mu(n)$ denotes the well known Möbius function (see \cite{A}), 
defined by $\mu(n)=(-1)^k$ if $n$ is a product of $k$ different prime 
numbers, and $\mu(n)=0$ if $n$ is not squarefree. 

We shall delete  the positivity condition, so
we shall say that the sequence of integers $(a_n)_{n=1}^\infty$ is  
\emph{pre-realizable} if  
$\sum_{d\mid n}\mu(n/d)a_d$ is divisible by $n$ for every natural 
number.

 In 1851 Kummer \cite{Km}  discovered what we call Kummer's 
 congruences for Bernoulli numbers, (see the book by Nielsen \cite{N}).
  Carlitz \cite{C} extended these 
 congruences to the generalized Bernoulli numbers of Leopoldt. Some 
 restrictions of Carlitz's results has been removed by the work of 
 Fresnel \cite{F}. These congruences are important for the definition 
 of the $p$-adic $L$-functions. 
 
 We establish a connection between these concepts that we can 
 formulate as in the following theorem.

\begin{theorem}
Let $(a_n)$ be a sequence that satisfies Kummer congruences for every 
rational prime, then for every natural number $b$ the sequence 
$(a_{b+n})_{n=1}^\infty$ is pre-realizable.
\end{theorem}

This theorem allow us to solve a problem posed by Gabcke \cite{G}. 
This is connected with the Riemann-Siegel formula.
In the investigation of the zeta function of Riemann it is important to 
compute the values of this function $\zeta(1/2+it)$ at points on the 
critical line with $t$ very high. Riemann found a very convenient formula 
for these computations, yet he does not publish anything about this 
formula. In 1932 C.~L.~Siegel was able to recover it from Riemann's nachlass. 
Now this formula is known as the Riemann-Siegel formula. 

To obtain the terms of this formula play a role certain numbers 
$\lambda_n$ that can be defined by a recurrence relation
\begin{equation*}
\begin{split}
\lambda_0&=1,\\
(n+1)\lambda_{n+1}&=\sum_{k=0}^n2^{4k+1}|E_{2k+2}|\lambda_{n-k}.
\end{split}
\end{equation*}
Here $E_{2n}$ denotes Euler numbers defined by 
\begin{equation*}
\frac{1}{\cosh x}=\sum_{n=0}^\infty \frac{E_n}{n!}x^n.
\end{equation*}
Hence $E_{2n+1}=0$ for $n\ge0$ and 
$$E_0=1,\quad E_2=-1,\quad E_4=5,\quad E_6=-61,\quad E_8=1\,385,\quad \dots$$

Gabcke  \cite{G} observed that the first six numbers $\lambda_n$ are 
integers and conjectured that this is so for all of them. 
Gabcke also considers  analogous sequences $(\varrho_n)$ and 
$(\mu_n)$. Although he does not mention it, the same motivations for 
his conjecture also supports that these too are integers sequences.
We prove all these conjectures. 
The proof of these assertions was the first motivation  of this  paper.

 In  \cite{PW} Puri and Ward ask if the sequence 
$(|E_{2n}|)_{n\ge1}$ is realizable. As we will see the solution of 
Gabcke's problem is connected with this one. We shall prove that in 
fact it is realizable.
\bigskip

Notations: When $p$ is a prime number and $m$ an integer we shall put 
$p^\alpha\parallel m$ to indicate that $p^\alpha$ is the greatest power 
of $p$ that divides $m$. We indicate this relation also by 
$\nu_p(m)=\alpha$.
We shall put $n\perp m$ to say that  $n$ and 
$m$ are relatively prime.


\section{Dynamical Zeta Function}

\begin{theorem}\label{th:beene}
    Let $(a_n)_{n=1}^\infty$ be a sequence of complex numbers and 
    define the sequence $(b_n)_{n=1}^\infty$ by
    \begin{equation}\label{eq:defb}
	nb_n=\sum_{d\mid n}\mu(n/d)a_d.
    \end{equation}
    Then  we have the equality between formal power series
    $$ \prod_{n=1}^\infty(1-x^n)^{-b_n}=
    \exp\Bigl(\sum_{n=1}^\infty\frac{a_n}n\Bigr).$$
\end{theorem}

\begin{proof}
    By the well known Möbius inversion formula the relation 
    (\ref{eq:defb}) is 
    equivalent to
    \begin{equation}\label{E:Mobius}
	a_n=\sum_{d\mid n} d b_d,
	\end{equation} therefore
we have the following equalities between formal power series
\begin{eqnarray}
\log\prod_{n=1}^\infty(1-x^n)^{-b_n}&=&-\sum_{n=1}^\infty b_n\log(1-x^n)=
\sum_{n=1}^\infty\sum_{k=1}^\infty b_n\frac{x^{nk}}{k}\nonumber\\
&=&
\sum_{m=1}^\infty\frac{x^m}{m} \Bigl(\sum_{n\mid m}nb_n\Bigr)=
\sum_{m=1}^\infty\frac{x^m}{m}  a_m.\nonumber
\end{eqnarray}
And this is equivalent to the equality we want to prove.
\end{proof}

\begin{theorem}\label{th:Aene}
     Let $(a_n)_{n=1}^\infty$ be a sequence of complex numbers and 
    define the sequence $(A_n)_{n=0}^\infty$ by the recurrence 
    relation
    \begin{equation*}\begin{split}
        A_{0}&=1,\\
        (n+1)A_{n+1}&= \sum_{k=0}^n A_{n-k} a_{k+1},\qquad n\ge0.
    \end{split}\end{equation*}
    Then we have the equality between formal power series
    $$\sum_{n=0}^\infty A_n x^n=
     \exp\Bigl(\sum_{n=1}^\infty\frac{a_n}{n}x^n\Bigr).$$
\end{theorem}

\begin{proof}
    First we have the equality between formal power series
    $$\sum_{n=1}^\infty n A_n x^{n-1}=\Bigl(\sum_{n=0}^\infty A_n x^n\Bigr)
    \Bigl(\sum_{n=1}^\infty a_n x^{n-1}\Bigr)$$
    because by the hypothesis the coefficient of $x^n$ is equal in both members.
    
    Since $A_0=1$ integrating formally give us
    $$\log\Bigl(\sum_{n=0}^\infty A_n x^n\Bigr)=\sum_{n=1}^\infty \frac{a_n}{n}
    x^n.$$
    That is equivalent to the relation we wanted to prove.
\end{proof}

The following theorem gives various equivalent conditions for a 
sequence of integers
 $(a_n)_{n=1}^\infty$ to be pre-realizable.

\begin{theorem} \label{T:equiv}
Given a sequence $(a_n)_{n\ge1}$ of integers, the following 
conditions are equivalent:
\begin{enumerate}
\item[(a)] The numbers $(b_n)_{n\ge 1}$ defined by 
$$n b_n=\sum_{d\mid n} \mu(n/d) a_d$$
are integers for every $n\in\N$. 

\item[(b)]
The numbers $(A_n)_{n\ge0}$ defined by 
\begin{equation}\label{eq:defA}
\begin{split}
A_{0}&=1,\\
(n+1)A_{n+1}&= \sum_{k=0}^n A_{n-k} a_{k+1},\qquad n\ge0
\end{split}\end{equation}
are integers for every $n\ge0$.

\item[(c)]  For  every prime number $p$ and natural numbers $n$, $\alpha$ 
with $p\perp n$
we have 
$$a_{np^\alpha}\equiv a_{n p^{\alpha-1}}\pmod{p^\alpha}.$$
\end{enumerate}
\end{theorem}

\begin{proof}
First we prove the equivalence of (a) and (b). 

(a)$\Longrightarrow$(b).
Assume  (a). By the definition of the $(b_n)$ and
 Theorem \ref{th:beene} we have
$$\prod_{n=1}^\infty(1-x^n)^{-b_n}=
    \exp\Bigl(\sum_{n=1}^\infty\frac{a_n}n x^n\Bigr).$$
and by the condition (a) of the theorem the $b_n$ are integers.
Let $(A_n)_{n=0}^\infty$ the numbers defined  by \ref{eq:defA}, we 
have to show that they are integers.
By  Theorem \ref{th:Aene} these numbers satisfies the relation
    $$\sum_{n=0}^\infty A_n x^n=
     \exp\Bigl(\sum_{n=1}^\infty\frac{a_n}n x^n\Bigr).$$
Thus we have
$$\sum_{n=0}^\infty A_n x^n=\prod_{n=1}^\infty(1-x^n)^{-b_n}.$$
Expanding this product, since the $b_n$ are integers, we get that 
the $A_n$  are also integers. 
Hence we have that (a) implies (b).

(b)$\Longrightarrow$(a).
Now, by hypothesis the numbers $(A_n)_{n\ge0}$ are integers. We can 
determine inductively a unique sequence of integers $c_n$ such that 
$$\sum_{n=0}^\infty A_n x^n=\prod_{n=1}^\infty (1-x^n)^{-c_n}.$$

In the first step observe that the coefficients of $x$ in both members 
must be the same, hence 
$$A_1=c_1. $$

Then observe that 
$$(1-x)^{c_1} \Bigl(\sum_{n=0}^\infty A_n x_n\Bigr)=1+ \sum_{n=2}^\infty 
A^{(2)}_n x^n,$$
where the numbers $A^{(2)}_n$ are integers. 

Assume by induction that we have determined integers $c_j$, for $j=1$, $2$, \dots 
$n-1$ such that 
$$\prod_{j=1}^{n-1}(1-x^j)^{c_j}\Bigl(\sum_{n=0}^\infty A_n x^n\Bigr)=
1+\sum_{k=n}^\infty A^{(n)}_k x^k.$$
Then the $A^{(n)}_k$ are integers and we can define 
$c_n=A_n^{(n)}$, that satisfies the induction hypothesis. 
Now we have
$$\prod_{j=1}^\infty (1-x^j)^{c_j}\Bigl(\sum_{n=0}^\infty A_n 
x^n\Bigr)=1.$$
By the hypothesis and Theorem \ref{th:Aene}
$$\sum_{n=0}^\infty A_n x^n=
\exp\left(\sum_{n=1}^\infty\frac{a_n}{n} 
x^n\right).$$
Therefore
$$\prod_{n=1}^\infty (1-x^n)^{-c_n}=\sum_{n=0}^\infty A_n x^n=
\exp\left(\sum_{n=1}^\infty\frac{a_n}{n} 
x^n\right).$$
Now take logarithms in both members to obtain
$$\sum_{n=1}^\infty c_n\sum_{k=1}^\infty \frac{x^{kn}}k=
\sum_{n=1}^\infty\frac{a_n}{n} x^n.$$
Reasoning as in the proof of Theorem \ref{th:beene} we get
$$a_m=\sum_{n\mid m} n c_n.$$
Therefore by the Möbius inversion formula  $c_n=b_n$ 
the numbers  defined on condition (a), and by construction 
these numbers $c_n$ are integers. 
Thus we have proved (a).

(a)$\Longrightarrow$(c).
We know that condition (a) is equivalent to the existence of integers 
$b_n$ that satisfy the equation (\ref{E:Mobius}). 

Assume that $p$ is a prime number and $n$ and $\alpha$ natural numbers such 
that $p\perp n$. 
Then 
$$a_{np^{\alpha}}=\sum_{d\mid np^\alpha} d b_d=
\sum_{k=0}^\alpha\sum_{d\mid n} d p^k b_{d p^k}.$$
Analogously 
$$a_{np^{\alpha-1}}=\sum_{k=0}^{\alpha-1}\sum_{d\mid n} d p^k b_{d p^k}.$$
Therefore 
$$a_{np^{\alpha}}-a_{np^{\alpha-1}}=\sum_{d\mid n} d p^\alpha b_{d p^\alpha}
\equiv0\pmod{p^{\alpha}}.$$

(c)$\Longrightarrow$(a).
Let $n$ be an integer. We have to show that 
$$\sum_{d\mid n} \mu(n/d) a_d$$
is divisible by $n$. Let $p^\alpha\parallel n$, with $\alpha\ge1$,
 then $n=p^\alpha m$ with 
$p\perp m$. 

Since $\mu(k)\ne0$ only when $k$ is squarefree, we get 
\begin{multline*}
\sum_{d\mid n} \mu(n/d) a_d=\sum_{d\mid m} \mu(m/d) a_{d p^\alpha }-
\sum_{d\mid m}\mu(m/d) a_{d p^{\alpha-1}}\\=
\sum_{d\mid m}\mu(m/d)\bigl( a_{d p^\alpha } -a_{d p^{\alpha-1}}\bigr)
\equiv 0\pmod{p^{\alpha}}.
\end{multline*}
The sum is divisible for every primary divisor  of $n$, and therefore
divisible by $n$.
\end{proof}

\section{Kummer congruences}

In 1851 Kummer \cite{Km} proved  the following theorem:

\begin{theorem}[Kummer]
    Let   $p$ be a prime number. Assume that
    $$\sum_{n=0}^\infty  a_n\frac{x^n}{n!}=\sum_{k=0}^\infty 
    c_k(e^{bx}-e^{ax})^k,$$
    where $a$, $b$ and the $c_k$ are integral $\pmod p$. Then the $a_n$ 
    are integers 
    $\pmod p$ and for $e\ge 1$, $n\ge1$,  $m\ge0$, and 
    $p^{e-1}(p-1)\mid w$ we have
    \begin{equation}\label{E:Kummer}
	\sum_{s=0}^n(-1)^s\binom{n}{s}a_{m+sw}\equiv0
	\pmod{(p^m,p^{ne})}.
    \end{equation}
\end{theorem}

The congruences (\ref{E:Kummer}) are usually called Kummer congruences. 
We shall say that the sequence $(a_n)$ satisfies Kummer congruences if 
we have  (\ref{E:Kummer}) for every prime number $p$. By Kummer 
theorem these sequences exist, but we are interested in some 
particular sequences. 

\begin{theorem}\label{T:Euler}
    The sequence $(E_{2n})_{n=1}^\infty$ satisfies Kummer congruences.
    \end{theorem}
    
\begin{proof}
    Since
    $$\frac1{\cosh 
    x}=\frac2{e^x+e^{-x}}=\frac{2}{2+(e^{x/2}-e^{-x/2})^2},$$
    Kummer theorem proves that $(E_n)_{n=1}^{\infty}$ satisfies 
    Kummer congruence (\ref{E:Kummer}) for every odd prime number $p$. 
    Therefore the sequence $(E_{2n})_{n=1}^\infty$ satisfies 
    these congruences for every odd prime number $p$. This reasoning 
    can be found in Kummer \cite{Km}.
    
    The above procedure does not give the case $p=2$, but 
    Fresnel  \cite{F} has extended Kummer congruences for Euler 
    number, as we see in the following lines.
    
    Let $\chi$  the function from $\Z$ to $\{-1,0,1\}$, such that 
    $\chi(n)=0$ if $n$ is even, $\chi(n)=1$ if $n\equiv1\pmod4$ and 
    $\chi(n)=-1$ if $n\equiv3\pmod4$, then the generalized Bernoulli 
    numbers associated to this character, ---see \cite{F} for 
    details---
    are related to Euler numbers as 
    $$\frac{B^n(\chi)}{n}=-\frac{E_{n-1}}{2}.$$
    In \cite{F} p.~319 we found that, when $2^e\parallel w$, with $e\ge1$
    $$\sum_{s=0}^n(-1)^s\binom{n}{s}\frac{B_{m+sw}(\chi)}{m+sw}\equiv0
    \pmod{(2^{n(e+2)}, 2^{m-1})}.$$
    With a change of notation this is equivalent to 
    $$\sum_{s=0}^n(-1)^s\binom{n}{s}\frac{E_{2m+sw}}{2}\equiv0
    \pmod{(2^{n(e+2)}, 2^{2m})}.$$
    Obviously this implies that for $2^{e-1}\mid w$ we have
    $$\sum_{s=0}^n(-1)^s\binom{n}{s}E_{2(m+sw)}\equiv0
    \pmod{(2^{ne},2^m)}.$$
\end{proof}

The above theorem is a model of many more interesting examples.  In 
Carlitz \cite{C}, it is proved that if $\chi$ is a primitive 
character $\bmod f$, and $f$ is divisible by at least two distinct 
rational primes, then $B^n(\chi)/n$ is an algebraic integer and
$$\sum_{s=0}^n(-1)^s\binom{n}{s}\frac{B^{n+1+sw}(\chi)}{n+1+sw}\equiv0
\pmod{(p^n,p^{e n})},$$
if $p^{e-1}(p-1)\mid w$.

Thus the sequence $(a_n)_{n=1}^\infty$, with 
$a_n=B^{n+1}(\chi)/(n+1)$ satisfies Kummer congruences, if the 
character $\chi$ is real. When $\chi$ is complex the sequence  defined 
by $a_n=\Tr(B^{n+1}(\chi)/(n+1))$ satisfies Kummer congruences.

The sequences that satisfies Kummer congruences are 
pre-realizable, as we will see in 
the following theorem. 

\begin{theorem}\label{T:KummerK}
    Let $(a_n)_{n=1}^\infty $ a sequence that satisfies Kummer 
    congruences. Then 
    $$a_{b+np^\alpha}\equiv a_{b+np^{\alpha-1}}\pmod{p^\alpha},$$
    for every natural numbers $b$, $n$, $\alpha$ and prime number $p$ 
    such that $p\perp n$.
    
    That is to say that if $(a_n)$ satisfies Kummer congruences then  
    for every natural number $b$,
    the sequence
    $(a_{b+n})_{n=1}^\infty$ is pre-realizable.
    \end{theorem}
    
\begin{proof}
    By (\ref{E:Kummer}), with $n=1$ we have
    $$a_{m+p^{e-1}(p-1)}\equiv a_m\pmod{(p^m, p^e)}.$$
    Therefore, for every natural number $k$, and assuming $m\ge e$
    $$a_{m+kp^{e-1}(p-1)}\equiv a_m\pmod{ p^e}.$$
    
    Now take $m=b+np^{\alpha-1}$, $k=n$ and $e=\alpha$. If 
    $b+np^{\alpha-1}\ge \alpha$, we get
    $$a_{b+np^{\alpha-1}+np^{\alpha-1}(p-1)}\equiv 
    a_{b+np^{\alpha-1}}\pmod{p^\alpha}.$$
    Since 
    $p^{\alpha-1}\ge\alpha$ for $p$ prime and $\alpha\ge1$, the 
    condition is satisfied and we get 
    $$a_{b+np^\alpha}\equiv 
    a_{b+np^{\alpha-1}}\pmod{p^\alpha}.$$
\end{proof}

\section{Euler numbers as numbers of fixed points}

We are now in position to solve the problem posed by Puri and Ward in 
\cite{PW}, 
they ask if the sequence $(|E_{2n}|)_{n=1}^\infty$ is realizable.
We shall show that this is true.

\begin{theorem}\label{th:Ereal}
    There exists a map $T\colon X\to X$, such that 
$$|E_{2n}|=\Fix T^n.$$
\end{theorem}

\begin{proof}
    First we show that $|E_{2n}|$ is a pre-realizable sequence.
    By Theorem \ref{T:Euler} the sequence $(E_{2n})_{n=1}^\infty$ 
    satisfies Kummer congruences. Thus by Theorem \ref{T:KummerK}, 
    for $p\perp m$, 
    $$E_{2m p^{\alpha}}\equiv E_{2m p^{\alpha-1}}\pmod{p^\alpha}.$$
    Therefore,
    $$|E_{2m p^{\alpha}}|\equiv |E_{2m p^{\alpha-1}}|\pmod{p^\alpha}.$$
By  Theorem \ref{T:equiv}  it follows that the numbers $b_n$, defined 
by 
$$n b_n=\sum_{d\mid n}\mu(n/d)|E_{2d}|,$$
are integers.

Now we must show that the numbers  $b_n$ are non negative.
To this end we observe that
$$nb_n\ge |E_{2n}|-\sum_{d=1}^{n/2}|E_{2d}|.$$
Now we apply the well known formula
$$1\le 
\Bigl(\frac{\pi}{2}\Bigr)^{2d+1}\frac{|E_{2d}|}{(2d)!}=2\sum_{k=0}^\infty
\frac{(-1)^k}{(2k+1)^{2d+1}}\le 2.$$
Thus 
\begin{eqnarray}
n b_n&\ge& (2n)!\Bigl(\frac{2}{\pi}\Bigr)^{2n+1}-2\sum_{d=1}^{n/2}
(2d)!\Bigl(\frac{2}{\pi}\Bigr)^{2d+1}\nonumber\\
&\ge&(2n)!\Bigl\{
\Bigl(\frac{2}{\pi}\Bigr)^{2n+1}-2\frac{(n)!}{(2n)!}
\sum_{d=1}^\infty\Bigl(\frac{2}{\pi}\Bigr)^{2d+1}\Bigr\}.\nonumber
\end{eqnarray}
We can compute the last sum and we get $0.433\dots$, therefore
$$nb_n\ge (2n)!\Bigl\{
\Bigl(\frac{2}{\pi}\Bigr)^{2n+1}-\frac{(n)!}{(2n)!}\Bigr\}.$$
This is positive for $n\ge2$, and we have $b_1=1\ge0$.
\end{proof}

The first values of these three sequences in this case are the 
following:
\begin{center}
\begin{tabular}[c]{|c|c|c|c|c|c|c|c|c|c|}
\hline
$a_n$ &  & 1 & 5 & 61 & 1385 & 50521 & 2702765 & 199360981 & 
\dots \\
\hline
$b_n$ &  & 1 & 2 & 20 & 345 & 10104 & 450450 & 28480140 &    \dots \\
\hline
$A_n$ & 1 & 3 & 23 & 371 & 10515 & 461869 & 28969177 & 2454072147 & 
 \dots \\
\hline
\end{tabular}
\end{center}

\section{Solution of Gabcke problem}

By Theorem \ref{T:equiv} the assertion of Gabcke, 
---that the numbers $\lambda_n$ are integers---, is 
equivalent to say  that the
sequence $a_{k}=2^{4k-3}|E_{2k}|$ is pre-realizable. We shall show 
that in fact it is realizable.

If $(a_n)$ and $(a'_n)$ are realizable, then the 
sequence $(a_na'_n)$ is also realizable. In fact given $T\colon X\to X$  
and $T'\colon Y\to Y$ such that $a_n=\Fix T^n$ and $a'_n=\Fix T'{}^n$, 
then it is easy to see that $T\times T'\colon X\times Y\to X\times Y$ 
satisfies $a_na'_n=\Fix(T\times T')$. 

Therefore by Theorem \ref{th:Ereal}  what we need is to prove that the sequence 
$2^{4n-3}$ is realizable.  This follows from the following theorem.

\begin{theorem}
    Let $a$ and $b\in\N$ such that $b\mid a$ and for  every prime number $p\mid a$
    $p\mid (a/b)$. Then the sequence $a^n/b$ is realizable.
\end{theorem}

\begin{proof}
    By the result of Puri and Ward we must show that the sequence $a^n/b$ is 
    pre-realizable and that the corresponding $b_n$ are non-negative 
    integers.
    
    First let $p\perp a$ be a prime number, and let $n\perp p$ and 
    $\alpha$ be natural numbers. We must show that 
    $$a^{np^\alpha}/b\equiv a^{np^{\alpha-1}}/b\pmod {p^\alpha}.$$
    Since $b\perp p$, this is equivalent to 
    $$a^{np^\alpha}\equiv a^{np^{\alpha-1}}\pmod {p^\alpha}.$$
    Now for $\alpha=1$ this is Fermat's little theorem, and for a 
    general $\alpha$ it follows, by induction, from the fact that
    for $\alpha\ge1$ if $a\equiv b\bmod {p^{\alpha}}$, then
    $a^p\equiv b^p\bmod{p^{\alpha+1}}$.
    
    Now if $p\mid a$, assume that $p^r\parallel a$ and $p^s\parallel b$. 
    By hypothesis we have $r\ge s+1$.
    We have to show that 
    $$a^{np^\alpha}/b\equiv a^{np^{\alpha-1}}/b\pmod {p^\alpha},$$
    where $p\perp n$ and $\alpha\ge1$. But the two numbers are 
    divisible by $p^{rnp^{\alpha-1}-s}$. All we have to show is that 
    $rnp^{\alpha-1}\ge s+ \alpha$. We can assume that $n=1$. 
    For $\alpha=1$ this is  $r\ge 
    s+1$ that is true by hypothesis. For other values of $\alpha$, 
    $\alpha\ge2$ and we 
    have
    $$rp^{\alpha-1}=r(p^{\alpha-1}-1)+r\ge (\alpha-1)+(s+1).$$
    
    Now we define the numbers $b_n$ by 
    $$nb_n=\sum_{d\mid n} \mu(n/d) a^d/b.$$
    By the previous reasoning we know the $b_n$ are integers.
    If $a=1$, it is easy to see that $b_1=1$ and $b_n=0$ for $n>1$. In 
    other case $a\ge2$ and
    we have 
    $$nb_n\ge \frac1b \Bigl( a^n-\sum_{d=1}^{n/2} a^d\Bigr).$$
    This is easyly seen to be non-negative.
\end{proof}

\begin{corollary}\label{th:cor}
The sequence $(a_n)_{n=1}^\infty$, where $a_n=2^{4n-3}$ is realizable.
\end{corollary}

The three sequences associated to this realizable sequence are 
\begin{center}
\begin{tabular}[c]{|c|c|c|c|c|c|c|c|c|c|c|}
\hline
$a_n$ &  & 2 & 32 & 512 & 8192 & 131072 & 2097152 & 33554432 & 
536870912 &  \dots \\
\hline
$b_n$ &  & 2 & 15 & 170 & 2040 & 26214 & 349435 & 4793490 & 67107840 &  \dots \\
\hline
$A_n$ & 1 & 2 & 18 & 204 & 2550 & 33660 & 460020 & 6440280 & 
91773990 & \dots \\
\hline
\end{tabular}
\end{center}

\medskip

Now we are in position to prove Gabcke's conjecture.
\begin{theorem}
    Let $\lambda_n$ the numbers defined by
    \begin{equation}\label{E:deflambda}
    \begin{split}
    \lambda_0&=1,\\
		(n+1)\lambda_n&= \sum_{k=0}^n 2^{4k+1}|E_{2k+2}|\lambda_{n-k},
		\qquad (n\ge0),
    \end{split}
    \end{equation}
$\varrho_n$ those defined by 
    \begin{equation}\label{E:defvarrho}
    \begin{split}
    \varrho_0&=-1,\\
    (n+1)\varrho_n&= -\sum_{k=0}^n 2^{4k+1}|E_{2k+2}|\varrho_{n-k},
    \qquad (n\ge0),
    \end{split}
    \end{equation}
and finally let $\mu_n=(\lambda_n+\varrho_n)/2$. All those numbers are 
integers.
    \end{theorem}
    
\begin{proof}
    By Theorem \ref{th:Ereal}, and Corollary \ref{th:cor} the 
    sequences $(|E_{2n}|)_{n=1}^\infty$ and $(2^{4n-3})_{n=1}^\infty$ 
    are realizable. Since the product of two realizable sequences is 
    realizable, the sequence
    $(2^{4n-3}|E_{2n}|)_{n=1}^\infty$ is realizable. Therefore it 
    satisfies condition (a) of Theorem \ref{T:equiv}. So it satisfies 
    condition (b), but this is precisely that the numbers $\lambda_n$ 
    are integers.
    
    Now condition (c) of the same Theorem gives us that with 
    $a_n=2^{4n-3}|E_{2n}|$  we have for every prime number $p$ and 
    natural numbers $n\perp p$ and $\alpha$ that
    $$a_{np^\alpha}\equiv a_{n p^{\alpha-1}}\pmod {p^\alpha}.$$
    Thus the same congruences are satisfied by the numbers $a'_n=-a_n$.
    Once again Theorem \ref{th:Ereal} says that the numbers $a'_n$ 
    satisfies condition (b). This is the same as saying that the 
    numbers $A'_n$ defined by 
    
        \begin{minipage}{1cm}\begin{eqnarray}\end{eqnarray}\end{minipage}\hfill
    \begin{minipage}{10cm}
	    \begin{eqnarray*}
		A'_0&=&1\\
		(n+1)A'_n&=& -\sum_{k=0}^n 2^{4k+1}|E_{2k+2}|A'_{n-k}
		\qquad (n\ge0)
	    \end{eqnarray*}\end{minipage}

	    \noindent are integers. But it is easyly seen that 
	    $\varrho_n=-A'_n$.

The affirmation about the numbers $\mu_n$ follows from the fact that 
$\lambda_n\equiv\varrho_n\pmod{2}$. That we prove in  Theorem  \ref{T:oh!}.
	    \end{proof}
	    
	    The following theorem is well known. I give a proof for 
	    completeness. 

\begin{theorem}
Let $s(n)$ be the sum of the digits of the binary representation of 
$n$, then
$$s(n)=n-\sum_{j=1}^\infty \Bigl\lfloor\frac{n}{2^j}\Bigr\rfloor.$$
\end{theorem}

\begin{proof}
Let the binary representation of $n$ be of type $\cdots 
0\overbrace{11\cdots1}^{\text{$k$ times}}$, with $k\ge0$, then 
$n+1=\cdots 
1\overbrace{00\cdots0}^{\text{$k$ times}}$. Therefore
$$s(n)-k=s(n+1)-1.$$
Also $k=\nu_2(n+1)$ the exponent of $2$ in the prime factorization 
of $n+1$.  

Thus we have proved that for every integer  $n\ge0$ 
\begin{equation}\label{E:referencia}
s(n+1)+\nu_2(n+1)=s(n)+1.
\end{equation}
We add this equalities for $n=0$, $1$, \dots, $n-1$ to get 
$$s(n)+\sum_{k=1}^n \nu_2(k)=n.$$

It is easily checked that 
$$\sum_{k=1}^n \nu_2(k)=\sum_{j=1}^\infty 
\Bigl\lfloor\frac{n}{2^j}\Bigr\rfloor.$$
\end{proof}

\begin{theorem} \label{T:oh!}
The numbers $\lambda_n$ and $\varrho_n$ defined by Equations (\ref{E:deflambda})  
and (\ref{E:defvarrho})
satisfy 
$$\nu_2(\lambda_n)=\nu_2(\varrho_n)=s(n).$$
\end{theorem}

\begin{proof}
First consider the sequence $\lambda_n$. Clearly the theorem is true 
for the first $\lambda_n$ which are
$$\lambda_0=1,\quad \lambda_1=2, \quad \lambda_3=82,\quad \lambda_4= 
10572.$$
Since Euler numbers $E_{2k}$ are odd, from the definition of 
$\lambda_n$ it follows that 
\begin{equation}\label{E:induccion}
\nu_2(n+1)+\nu_2(\lambda_{n+1})=\nu_2\Bigl(
\sum_{k=0}^n 2^{4k+1}|E_{2k+2}|\lambda_{n-k}\Bigr).
\end{equation}
By induction  the terms of this sum are exactly divided by the powers 
of $2$ of exponents
$$1+s(n),\quad 5+ s(n-1),\quad 9+s(n-2),\quad \dots \quad (4n+1)+s(0)$$
This is a strictly increasing sequence, since 
$$s(n)-s(n-1)=1-\nu_2(n)<4.$$

Hence from (\ref{E:induccion}) we get
$$\nu_2(n+1)+\nu_2(\lambda_{n+1})=1+s(n).$$
By (\ref{E:referencia})
$$\nu_2(\lambda_{n+1})= s(n)-\nu_2(n+1)+1=s(n+1).$$

The same proof applies to the sequence $(\varrho_n)$.
\end{proof}

\section{Examples}

We give here some examples of numbers satisfying our Theorem \ref{T:equiv}. 
First consider the  case of the numbers of Gabcke $A_n=\lambda_n$. The 
first terms of the associated sequences are given by the following 
table.
\begin{center}
\begin{tabular}[c]{|c|c|c|c|c|c|c|c|c|}
\hline
$a_n$  & & 2 & 160 & 31232 & 11345920 & 947622146676 & 957663025230936 & \dots \\
\hline
$b_n$ &  & 2 & 79 & 10410 & 2836440 & 1324377702 & 944684832315 & \dots \\
\hline
$\lambda_n$ & 1 & 2 & 82 & 10572 & 2860662 & 1330910844 & 947622146676 & \dots \\
\hline
\end{tabular}
\end{center}
\medskip

We can give arbitrarily a sequence of integers $(b_n)$ and  obtain 
sequences $(a_n)$ and $(A_n)$ that automatically satisfy our 
theorems. We give two simple examples. 

With $b_n=1$ for every $n$,  we get  $a_n =\sigma(n)$.
\begin{center}
\begin{tabular}[c]{|c|c|c|c|c|c|c|c|c|c|c|c|c|c|c|}
\hline
$a_n$  & & 1 & 3 & 4 & 7 & 6 & 12 & 8 & 15 & 13 & 18 & 12 & 28 & \dots \\
\hline
$b_n$ &  & 1 & 1 & 1 & 1 & 1 & 1 & 1 & 1 & 1 & 1 & 1 & 1 &  \dots \\
\hline
$A_n$ & 1 & 1 & 2 & 3 & 5 & 7 & 11 & 15 & 22 & 30 & 42 & 56 & 77 & \dots \\
\hline
\end{tabular}
\end{center}

With $b_n=-24$, the numbers $A_n$ are given by   Ramanujan's 
$\tau$ function
 $A_n=\tau(n+1)$.
\begin{center}
\begin{tabular}[c]{|c|c|c|c|c|c|c|c|c|c|c|}
\hline
$a_n$  & & -24 & -72 & -96 & -168 & -144 & -288 & -192 & -360 &   \dots \\
\hline
$b_n$ & & -24 & -24 & -24 & -24 & -24 & -24 & -24 & -24 &    \dots \\
\hline
$A_n$ & 1 & -24 & 252 & -1472 & 4830 & -6048 & -16744 & 84480 & -113643 &  \dots \\
\hline
\end{tabular}
\end{center}
\medskip

Finally let $(T_n)$ be  the tangent numbers with the notation of 
\cite{K}. We have $T_{2n}=0$ 
$$T_{2n+1}=(-1)^{n} \frac {4^{n+1}(4^{n+1}-1)B_{2n+2}}{2n+2}.$$
It can be proved that $a_n=(-1)^nT_{2n+1}$ satisfies Kummer 
congruences. It follows that the sequence $(T_{2n+1})$ is realizable, in 
this case the three sequences are
\begin{center}
\begin{tabular}[c]{|c|c|c|c|c|c|c|c|c|c|}
\hline
$a_n$  & & 2 & 16 & 272 & 7936 & 353792 & 22368256 & 1903757312 &    \dots \\
\hline
$b_n$ & & 2 & 7 & 90 & 1980 & 70758 & 3727995 & 271965330 &     \dots \\
\hline
$A_n$ & 1 & 2 & 10 & 108 & 2214 & 75708 & 3895236 & 280356120 &   \dots \\
\hline
\end{tabular}
\end{center}

\bibliographystyle{amsplain}

\end{document}